\documentclass[letterpaper,11pt,oneside]{amsart}
\usepackage{amsmath}
\usepackage{amssymb}
\usepackage{latexsym}
\usepackage{amsfonts}
\usepackage[latin1]{inputenc}
\usepackage[english]{babel}
\usepackage{multicol}
\usepackage{oldgerm}
\usepackage{amscd,tikz}
\usepackage{enumerate}

\newtheorem{teo}{Theorem}

\begin{document}

\title{Quadratic forms representing the $p$th Fibonacci number}

\author{Pedro Berrizbeitia}
\address{Departamento de Matem\'aticas Pura y Aplicada\\
Universidad Sim\'on Bol\'{\i}var\\
Caracas, Venezuela}
\email{pberrizbeitia@gmail.com}

\author{Florian Luca}
\address{School of Mathematics\\
University of the Witwatersrand\\
P. O. Wits 2050, South Africa}
\email{florian.luca@wits.ac.za}

\author{Alberto Mendoza}
\address{Departamento de Matem\'aticas Pura y Aplicada\\
Universidad Sim\'on Bol\'{\i}var\\
Caracas, Venezuela}
\email{jacob@usb.ve}


\date{\today}

\begin{abstract}

In this paper, we show that if $p\equiv 1\pmod 4$ is prime, then $4F_p$ admits a representation of the form $u^2-pv^2$ for some integers $u$ and $v$, where 
$F_n$ is the $n$th Fibonacci number. We prove a similar result when $p\equiv -1\pmod 4$.

\medskip



\end{abstract}

\maketitle


\section{Introduction}

Let $\{F_n\}_{n\ge 0}$ be the Fibonacci sequence given by $F_0=0,~F_1=1$ and $F_{n+2}=F_{n+1}+F_n$ for all $n\ge 0$. 
The starting point for the investigation of the subject in the title is the formula 
\begin{equation}
\label{eq:first}
F_{2n+1}=F_n^2+F_{n+1}^2 
\end{equation}
known to Lucas (take $Q=-1$ in formula (34) in Lucas's seminal 1878 paper \cite{Lucas}) since it implies that every Fibonacci number of odd index can be represented as the sum of two squares of integers numbers. This is a question which leads naturally to the investigation of Fibonacci numbers $F_n$ which can be represented under the form $u^2+dv^2$ with some integers $u$ and $v$ and some integer $d$ which is either  fixed or depends on $n$. For example, in \cite{Savin}, it is shown that if $n\equiv 7\pmod {16}$, then $F_n=u^2+9v^2$ holds with some positive integers $u$ and $v$. For general results 
regarding this problem when $d$ is fixed, see \cite{BL}.

In \cite{AL}, it was noted that if $n=p^2$ is the square of an odd prime $p\ne 5$, then $p$ divides $F_{\frac{p^2-1}{2}}$, hence formula \eqref{eq:first}  implies that  $F_{p^2}= u^2 + p^2 v^2$, for some integers $u$ and $v$. Motivated by this  observation, they introduced and estimated the counting function of the infinite set 
$$
S= \{ n~:~F_n = u^2 + n v^2~{\text{\rm with~some~integers}}~ u,v\}.
$$ 
In the course of their investigation, they found computational evidence that indicated that every prime $p \equiv 1 \pmod 4$  belongs to $S$. In  \cite{ABL}, it was proved this fact is true; that is that  if $p\equiv 1\pmod 4$, then $F_p=u^2+pv^2$ for some integers $u$ and $v$. The proof makes use of basic facts in Galois Theory and basic properties of the norm function of finite extensions of ${\mathbb Q}$. 
Prior, it was shown in \cite{Savin} that the above formula never holds if instead of $p\equiv 1\pmod 4$, we have $p\equiv 3,7\pmod {20}$. In this paper, we find quadratic forms representing $F_p$, for all primes $p \equiv -1 \pmod 4$. In fact, we extend the method developed in \cite{ABL} to study other quadratic form representations for $F_p$ or $4 F_p$, all depending only on the congruence modulo $4$ of the prime $p$.
We state the main result of this paper:

\begin{teo}
\label{thm:main}
\begin{itemize}
\item[(i)] If $p\equiv 1\pmod 4$ is prime, then $4F_p=u^2-pv^2$ holds with some (in fact, infinitely many) pairs of positive integers $(u,v)$.
\item[(ii)] If $p\equiv -1\pmod 4$, then $4F_p=5u^2+pv^2$ holds with some integers $u$ and $v$.
\end{itemize}
\end{teo}
 
\section{The proof of Theorem \ref{thm:main}}

For a positive integer $n$ let $\zeta_n$ be a primitive $n$th root of unity. We let 
$$
(\alpha,\beta)=({(1+{\sqrt{5}})}/{2}, {(1-{\sqrt{5}})}/{2}),
$$ 
and use the fact that 
\begin{equation}
\label{eq:1}
F_n=\frac{\alpha^n-\beta^n}{\alpha-\beta}\quad {\text{\rm holds~for~all}}\quad n\ge 0.
\end{equation}

Our proof proceeds by noting that formula \eqref{eq:1} entails
\begin{equation}
\label{eq:F}
F_n=\prod_{t=1}^{n-1} (\alpha-\beta \zeta_n^t).
\end{equation}
We assume that $n=p$ is an odd prime and denote the field ${\mathbb Q}({\sqrt{5}},\zeta_p)$ by ${\mathbb L}.$ 
We let $p^*=(-1)^{\frac{p-1}{2}}p$. 
For any subfield $\mathbb{F}$ of $\mathbb{L}$, and for any $\gamma \in \mathbb{L}$, $N_{\mathbb{L}/\mathbb{F}}(\delta)$ will denote the relative norm of $\gamma$ from $\mathbb{L}$ to $\mathbb{F}$. So $N_{\mathbb{L}/\mathbb{F}}(\gamma)=\prod_{\sigma \in {\text{\rm Gal}}(\mathbb{L}/\mathbb{F})}  \sigma(\gamma),$ where ${\text{\rm Gal}}(\mathbb{L}/\mathbb{F})$ is the Galois Group of $\mathbb{L}$ over $\mathbb{F}$.

Note that $\mathbb{L}$ is a number field of  degree $2(p-1)$ over $\mathbb{Q}$, that is, $[{\mathbb L}:{\mathbb Q}]=2(p-1)$. 

Note further that 
$
{\text{\rm Gal}}({\mathbb L}/{\mathbb Q})=\langle \tau\rangle \times \langle \sigma_5\rangle,
$
where
\begin{itemize}
\item $\tau(\zeta_p)=\zeta_p^{g}$ for some generator $g$ of ${\mathbb Z}_p^*$, and  $\tau({\sqrt{5}})={\sqrt{5}}$;
\item $\sigma_5(\zeta_p)=\zeta_p$, and $\sigma_5({\sqrt{5}})=-{\sqrt{5}}$.
\end{itemize}
For integers $d$ and $d'$ we let ${\mathbb Q}_d={\mathbb Q}({\sqrt{d}})$ and we let ${\mathbb Q}_{d,d'}={\mathbb Q}({\sqrt{d}},{\sqrt{d'}} )$. 
The following diagram of subfields of $\mathbb{L}$  is useful:

\bigskip

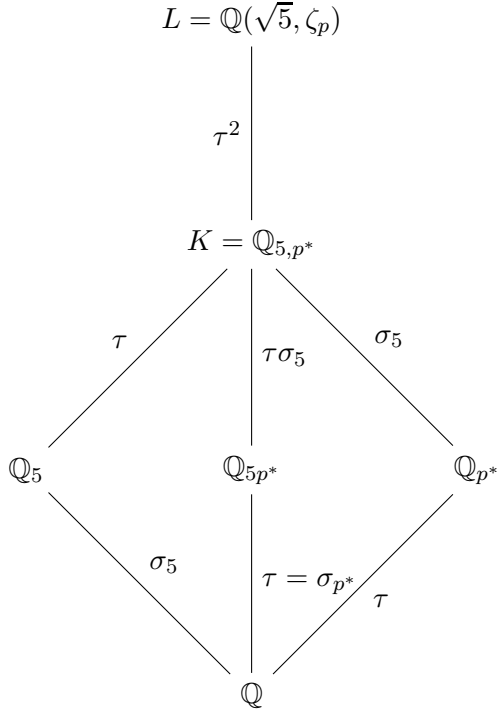
\begin{tikzpicture}[node distance=3cm, auto]
 \node (L) {$L=\mathbb{Q}(\sqrt{5},\zeta_p)$};
  \node (Q5p^*) [below of=L] {$K=\mathbb{Q}_{5,p^*}$};
  \node (Q5) [below of=Q5p^*, left of=Q5p^*] {$\mathbb{Q}_5$};
  \node (Q5p) [below of=Q5p^*] {$\mathbb{Q}_{5p^*}$};
  \node (Qp) [below of=Q5p^*, right of=Q5p^*] {$\mathbb{Q}_{p^*}$};
  \node (Q) [below of=Q5p] {$\mathbb{Q}$};
 \draw[-] (L) to node [swap] {$\tau^2$} (Q5p^*);
  \draw[-] (Q5p^*) to node [swap] {$\tau$} (Q5);
  \draw[-] (Q5p) to node {$\tau=\sigma_{p^*}$} (Q);
  \draw[-] (Q5p^*) to node {$\sigma_5$} (Qp);
  \draw[-] (Q5p^*) to node  {$\tau\sigma_5$} (Q5p);
  \draw[-] (Q5) to node {$\sigma_5$} (Q);
  \draw[-] (Qp) to node {$\tau$} (Q);
\end{tikzpicture}

\bigskip

From the diagram, one can see that
$$
{\text{\rm Gal}}({\mathbb {L}}/{\mathbb{K}})=\langle \tau^2\rangle,
$$ 
 $$
{\text{\rm Gal}}({\mathbb{L}}/{\mathbb {Q}_5})=\langle \tau\rangle,
$$ 
$$
{\text{\rm Gal}}({\mathbb{L}}/{\mathbb {Q}_{p^*}})=\langle \tau^2\rangle \times \langle \sigma_5 \rangle,
$$  
$$
{\text{\rm Gal}}({\mathbb{L}}/(\mathbb {Q}_{5p^*})=\langle \tau\sigma_5 \rangle.
$$ 
From (2), we obtain
$$F_p=\prod_{t=1}^{p-1} (\alpha-\beta \zeta_p^t)=\prod_{s=1}^{p-1} (\alpha-\beta \tau^s(\zeta_p))= N_{\mathbb{L}/\mathbb{Q}_5}(\alpha - \zeta_p \beta).$$
We let $\gamma = \alpha - \zeta_p \beta$, so $\gamma \in \mathbb{L}$.
We also let  
$$ 
\Gamma= N_{\mathbb{L}/\mathbb{K}}(\gamma) = \prod_{t=1}^{\frac{p-1}{2}} (\alpha-\tau^{2t} (\zeta_p)\beta)=\prod_{r \in R} (\alpha-\zeta^r \beta),
$$ 
where $R$ is the set of quadratic residues in $\mathbb{Z}_p^*.$
Again from (2), we get 
$$
F_p =\prod_{t=1}^{p-1} (\alpha- \zeta_p^t \beta )=\prod_{r \in R} (\alpha-\zeta_p^r \beta) \prod_{r \in R} (\alpha-\zeta_p^{gr} \beta),
$$ 
where $\tau(\zeta_p)=\zeta_p^g$. It follows that  
 $F_p=\Gamma \tau(\Gamma)$.
Next we compute $\sigma_5(\Gamma)$ and get: 
$$
\sigma_5(\Gamma)= \prod_{r \in R} \sigma_5(\alpha - \zeta_p^r \beta) = \prod_{r \in R} (\beta - \zeta_p^r \alpha) = \prod_{r \in R}\left (\frac{-1}{\alpha} + \frac{1}{\zeta_p^{-r} \beta}\right)= \frac{ \prod_{r \in R} (\alpha - \zeta_p^{-r} \beta)}
{(-1)^{\frac{p-1}{2}} \prod_{r \in R} \zeta_p^r},$$ 
where we have used for the last three equalities that
$$\sigma_5 (\alpha)=\beta,\quad \sigma_5 (\beta)=\alpha,\quad \sigma_5(\zeta_p)=\zeta_p\quad \alpha \beta = -1,
$$ which are all trivial to verify. 
This leads to the following:
\begin{enumerate}
\item if $p \equiv 1 \pmod 4$ then $\sigma_5(\Gamma) = \Gamma$;
\item if $p\equiv -1 \pmod 4$ then $\sigma_5(\Gamma) = -\tau(\Gamma)$.
\end{enumerate}

\subsection{Case i) of Theorem 1}

In that case $\sigma_5(\Gamma) = \Gamma$, so $\Gamma$ is fixed by $\tau^2$ and by $\sigma_5$. It follows that  
$\Gamma \in \mathbb{Q}_{p^*}$. Hence, $\Gamma = x + y \sqrt {p}$ for some half integers $x$ and $y$.
Since $F_p = \Gamma \tau(\Gamma)$, then $F_p = x^2 - p y^2$. Multiplying $\Gamma$ by norm 1 units of the ring of integers of $\mathbb{Q}_{p}$ leads to infinitely many solutions.

\subsection{Case ii) of Theorem 1}

Now $p^*=-p$ and  $\sigma_5(\Gamma) = -\tau(\Gamma)$. 
Since we have $\Gamma \in K=\mathbb{Q}_{5,-p},$ we get that 
\begin{equation}
\label{eq:Gamma1}
\Gamma = r + s\sqrt{5} + t{\sqrt {-p}} + u\sqrt{-5p}
\end{equation} 
with rational  numbers $r$, $s$,  $t$, $u$ whose denominator divides $4$. The condition 
$$
\sigma_5(\Gamma) = -\tau(\Gamma)  \quad {\text{\rm implies}}\qquad r = u = 0.
$$ 
It follows that  $\Gamma = s \sqrt 5 + t \sqrt -p$. Apriori, $4s$ and $4t$ are integers. We need to be more precise and justify that $2s$ and $2t$. In order to do so, we argue as follows. The quadratic fields ${\mathbb Q}_5$ and ${\mathbb Q}_{-p}$ have integral bases $\{1,(1+{\sqrt{5}})/2\}$, (respectively $\{1,(1+{\sqrt{-p}})/2\}$) and coprime discriminants $5$ and $-p$. 
It follows that an integral basis for ${\mathbb Q}_{5,-p}$ is 
$$
\{1,(1+{\sqrt{5}})/2, (1+{\sqrt{-p}})/2,(1+{\sqrt{5}})(1+{\sqrt{-p}})/4\}
$$
(see Exercise 4.5.13 in \cite{Mur}). Hence, for some integers $a,b,c,d$, we have
\begin{eqnarray}
\label{eq:Gamma2}
\Gamma & = & a+b(1+{\sqrt{5}})/2+c(1+{\sqrt{-p}})/2+d(1+{\sqrt{5}})(1+{\sqrt{-p}})/4\nonumber\\
& = & (a+b/2+c/2+d/4)+(b/2+d/4){\sqrt{5}}+(c/2+d/4){\sqrt{-p}}+(d/4){\sqrt{-5p}}.
\end{eqnarray}
Identifying coefficients in \eqref{eq:Gamma1} and \eqref{eq:Gamma2}, we get 
$$
r=a+b/2+c/2+d/4,\quad s=b/2+d/4,\quad t=c/2+d/4,\quad u=d/4.
$$
Since $r=u=0$, we get $d=0$, therefore $2s=b$ and $2u=c$ are integers. Since $F_p = \Gamma \tau(\Gamma)= 5 s^2 + p t^2$, we get the desired result.

\medskip

\end{document}